\documentclass{article}

\usepackage{amsfonts}
\usepackage{amsthm}
\usepackage{latexsym}

\newtheorem{thm}{Theorem}
\newtheorem{define}{Definition}
\newtheorem{lemma}{Lemma}

\newtheorem{fct}{Fact}

\theoremstyle{remark}
\newtheorem*{rem}{Remark}

\begin{document}

\title{ C-Groups\thanks{Research funded by grants from the 
NSF grant number NSF(DMS-0097113), and  NSA, grant number NSA(MSPR-01IC-180).}}
\author{Kazem Mahdavi\thanks{Project Mentor, mahdavk@potsdam.edu}, 
Margrethe Flanders\thanks{Group Leader, mebflanders@att.net}, 
Avi Silterra\footnote{ajs2@andrew.cmu.edu}, 
Mihai Tohaneanu\footnote{mihai@dartmouth.edu}}
\maketitle

\begin{abstract}
In this paper we explore the structure and 
properties of C-groups. We define
 a C-group as 
a group $G$ with $rk(G) < rk(Z(G))$.  Using GAP 
(a group theory program) 
and traditional methods, we identified 
an interesting infinite class of
 C-groups. In particular,
 we have proved that there is always a 
C-group of order an integer multiple 
of a fifth power 
of a prime. One way to obtain C-groups 
is take the direct product 
of certain C-groups, mentioned in this paper, 
with other appropriate groups.
  We, at this time, do not know whether the 
C-groups discussed in this paper 
are the building block of all finite C-groups.  
But a complete 
classification of finite or infinite 
C-groups is an interesting 
problem. We have also formulated a 
number of open questions relating to
 C-groups: Are they 
all solvable? What is the structure of 
the C-groups that are not in our
 class?  
Is the minimal number 
of generators of the center always polynomially
 bounded by the minimal number of 
generators of the
 group?  What are the isoperimetric 
inequalities of infinite C-groups?
\end{abstract}

\section*{Introduction}

Call the minimum cardinality of a generating set 
for a group $G$ the rank of $G$
and denote it by $rk(G)$.

In this paper we study a class of finitely generated groups $G$ such that 
$rk(G) < rk(Z(G))$ where $Z(G)$ is the center of $G$. We call these 
groups \textbf{C-groups}. 

We came across this class of groups while trying to study finitely
 generated groups with 
infinitely generated centers.  One class of such groups is the set 
of $4\times 4$ 
upper-triangular matrices over the set of rationals with a power of
 a prime $p$ in 
the denominator:
\[ \left( \begin{array}{cccc}
1 & * & * & * \\
0 & u_1 & * & * \\
0 & 0 & u_2 & * \\
0 & 0 & 0 & 1
\end{array} \right) \]
where $u_i$ denotes a positive unit. (For further details, see \cite{abe}).

We used GAP, a system for computational group theory, both to find finite
 C-groups of order 
less than 1000 (with the exceptions of 512 and 768) and to gain insight
 into their structures
 \cite{gap}. In this paper we explain the results of these investigations.

\section{Definitions and Notation}


\begin{define}[Elementary abelian]
A group $G$ is elementary abelian if $G = \prod_{i < n} \mathbb{Z}_p$ for
some $n \in \mathbb{N}$ and $p$ a prime.
\end{define}

\begin{define}[Frattini Subgroup]
The \textbf{Frattini subgroup} of $G$, denoted $\Phi (G)$, is the
 intersection of all 
maximal subgroups of G when they exist and $G$ otherwise.
\end{define}

\begin{define}[Nilpotent]
Let $G$ be a group, 
define $Z_1(G)=Z(G)$ and 
$Z_i(G)$ is the inverse image of $Z(G/Z_{i-1}(G))$ under the 
canonical projection 
$G\rightarrow G/Z_{i-1}(G)$. The group $G$ is \textbf{nilpotent} (of class $n$)
if $Z_n(G)=G$ for some $n$.
\end{define}

\begin{define}[Commutator]
Let $G$ be a group and let $x,\:y \in G$. Then 
$[x,y]=xyx^{-1}y^{-1}$ is called the 
\textbf{commutator} of $x$ and $y$. The subgroup of $G$ generated by the set 
$\{[x,y] \mid x,:y \in G \}$ is called the \textbf{commutator subgroup} of $G$ 
and will be denoted $G'$.
\end{define}

\begin{define}[Solvable]
Let $G$ be a group, and let $G^{(1)}:=G'$. Then for $i\geq 1$, define 
$G^{(i)}$ by $G^{(i)}=(G^{(i-1)})'$. The group $G$ is said to be
 \textbf{solvable} if $G^{(n)}=\langle e\rangle$ for some $n$.
\end{define}

\begin{define}[$p$-Group]
A group in which every element has order a 
power of some fixed prime $p$ is 
called a \textbf{$p$-group}.
\end{define}

\begin{lemma}
\label{one}
The set $\mathbb{Z}_{n_1} \times \mathbb{Z}_{n_2} \times \dots 
\times \mathbb{Z}_{n_k}$ 
with $\gcd(n_1,n_2,\ldots, n_k)>1$ has a minimum of $k$ generators.
\end{lemma}

\begin{lemma}
\label{two}
Any group of order $p^2$ for $p$ a prime is cyclic or elementary abelian.
\end{lemma}

\begin{thm}
A finite group $G$ is a $p$-group $\Leftrightarrow$ the 
order of $G$ is a power of some prime $p$.
\end{thm}

\begin{thm}
\label{nil}
Every finite $p$-group is nilpotent.
\end{thm}
For more information on nilpotent and solvable groups, 
see \cite{hun}, pages 100-107.

\begin{thm}
\label{frat}
\begin{description}
  \item[Burnside's Basis Theorem] \label{bbt} If $G$ is a finite p-group, 
    then $\forall g \in G, \ g^p \in \Phi(G)$, $G' \leq \Phi(G)$,
    $\bar{G} = G / \Phi(G)$ is elementary abelian, and 
    $rk(\bar{G}) = rk(G)$.
  \item[Universal Mapping Property of Frattini Subgroup] \label{umpfrat} If $H$ is a finite
	 p-group,
    $N \lhd H$ and $H / N$ is elementary abelian then $\Phi(H) \leq N$. 
\end{description}
\end{thm}
See \cite{dum}
or \cite{rob} for details.

\section{Main Results}

\begin{thm}[Main Result]
\label{main}
Let $G = \mathbb{Z}_{n_1} \times \mathbb{Z}_{n_2} \times \mathbb{Z}_{n_3}$ 
where $n_1 \mid n_2$, $n_1 \mid n_3$, and $\gcd( n_1, \frac{n_2}{n_1}, 
\frac{n_3}{n_1})~>~1$, 
and binary operation
\newline\indent\((x_1, y_1, z_1)\cdot(x_2, y_2, z_2) =\newline\indent\indent
( (x_1 + x_2 + y_2z_1)\bmod n_1, (y_1 + y_2)\bmod n_2, (z_1 + z_2)\bmod n_3 ).\)
\newline Then $(G,\cdot)$ is a C-group.
\end{thm}

\begin{proof} It is easy to verify that $G$ is a group under the above 
binary operation 
where $(0,0,0)$ is 
the identity element and \[(x,y,z)^{-1} = ((n_1-x) + yz, n_2-y, n_3-z).\]

The generators of $G$ are $a=(0,1,0)$ and $b=(0,0,1)$, since 
$(1,0,0)=bab^{-1}a^{-1}$,
 and 
they satisfy the 
following relations:
\[a^{n_2}=b^{n_3}=(0,0,0),\ (bab^{-1}a^{-1})^{n_1}=(0,0,0),
\ a^{n_1}b=ba^{n_1},
\ ab^{n_1}=b^{n_1}a.\]
Therefore, the canonical form of words in $G$ is 
\[a^{k_1n_1}b^{k_2n_1}[a,b]^{k_3}a^{k_4}b^{k_5}\]
where $k_1 = 0,1,\ldots,\frac{n_2}{n_1}$, $k_2 = 0,1,
\ldots,\frac{n_3}{n_1}$, and
 $k_3,k_4,k_5 = 0,1,...,n_1-1$. Thus the order of $G$ is $n_1n_2n_3$.

Let $C=Z(G)$. Then 
\[C=\{ (x,n_1y,n_1z) \mid x \in \mathbf{Z}_{n_1}, y \in \mathbf{Z}_{n_2}, 
z \in \mathbf{Z}_{n_3} \}.\] So, 
$C \cong \mathbb{Z}_{n_1} \times 
\mathbb{Z}_\frac{n_2}{n_1} \times \mathbb{Z}_\frac{n_3}{n_1}$ 
of order $n_1(\frac{n_2}{n_1})(\frac{n_3}{n_1})$ by Lemma \ref{one}.

Therefore, $G$ is a C-group of order $n_1n_2n_3$.
\end{proof}

\begin{rem}
This defines an infinite class of C-groups which we call $\alpha$C-groups.
\end{rem}

\begin{thm}
There is a unique C-group of order $p^5$ for any prime $p$.
\end{thm}

\begin{proof}
Existence follows from Theorem \ref{main} with $n_1=p,\ n_2=p^2,\ n_3=p^2$.

Uniqueness is more difficult. Suppose $G$ is a C-group such that $|G|=p^5$.
First we'll show $rk(Z(G))=3$ and $rk(G)=2$. 
Let $|G/Z(G)|= p^q$. $q=1$ implies $G$ 
abelian, and $q \geq 3$ implies that $rk(Z(G))\leq 2$. 
Hence, $q=2$ and so $|Z(G)| = p^3$. Therefore, 
 $rk(Z(G)) \leq 3$, so as $G$ is a C-group 
we have $rk(Z(G))=3$ and $rk(G)=2$.

Next we will show that $\Phi(G) = Z(G)$, though we'll only need 
$\Phi(G) \leq Z(G)$. 
To begin, note that $G/Z(G)$ is elementary abelian by Lemma \ref{two}.
  By the UMP of the Frattini
subgroup, we have $\Phi(G) \leq Z(G)$. Next by the 
Burnside's Basis Theorem and $rk(G)=2$, we must have $|\Phi(G)|= p^3$ 
in order that $|G / \Phi(G)| = p^2$. Hence, $Z(G)=\Phi(G)$.

Choose $a,b \in G$ such that
$\langle a,b \rangle = G$. By the previous paragraph, 
$a^p, b^p \in Z(G)$ so $|a|, |b| \leq p^2$ since 
$Z(G)$ is elementary abelian.

To finish, 
$[a,b] \in Z(G)$ and it is easy to see 
any $x \in G$ has a representation $x = a^i b^j [a,b]^k$ for 
$i \leq |a|$, 
$j \leq |b|$, $k \leq |[a,b]|$. This is because $a, b$ generate
$G$, any instance of $aba$ in the representation for an element can
be replaced with $aab[a,b]$, and $[a,b] \in Z(G)$. 
Since $|G|=p^5$, and 
$|[a,b]|\leq p$ (as $Z(G)$ is elementary abelian), then $|a|=|b|=p^2$
and $|[a,b]|=p$. \\
Therefore, $G$ has the same presentation as our 
C-group of the same size, 
so the two are isomorphic.
\end{proof}

\begin{fct}
All C-groups of the class defined above are nilpotent of class 2.
\end{fct}

\begin{proof} Let $x=(x_1,x_2,x_3)$ and $y=(y_1,y_2,y_3)$ in some
 C-group $G$ of the class 
defined above. Then $[x,y]=(y_2x_3+(n_2-1)x_2y_3,0,0)$. So $[x,y]$ 
is clearly in the center 
of $G$. Therefore, all C-groups in our class are nilpotent of class 2.
\end{proof}



\section*{More C-Groups}
The multiplication defined for $\alpha$C-groups appears to be a semidirect product,
but GAP tests indicated that was not the case.
Though the number of $\alpha$C-groups is 
infinite and can be extended infinitely by 
right direct product with an appropriate abelian group such 
as $\mathbf{Z}_{kn_1}$ or 
$\mathbf{Z}$ where $k$ is a positive integer, 
this class does not describe all
C-groups.
For example, there is a C-group of order 64 that is nilpotent of class 3 with 
the following 
group representation 
\[ <a,b \mid a^4=b^8=1,\ a^2b=ba^2,\ b^2a=ab^{-2}, (b^{-1}a)^2=(ab)^2>. \]

There are also two C-groups of order 96
that are not nilpotent. These facts
 raise some 
interesting questions: Is there a bound on the 
nilpotency class of C-groups? What is the 
structure of the nilpotent C-groups not in our class?
 What is the structure of non-nilpotent
 C-groups? In addition, all C-groups of order 
less than 1000 are solvable. Are all C-groups 
solvable?

A brief examination of the GAP data reveals that the number of 
generators of the center of 
C-groups of order less than 1000 is a linear function of the 
number of group generators. 
We conjecture that for finite C-groups, the number of generators of the center is
 polynomially
 bounded by the number of generators of the group. 
Clearly, there is no such bound for
 C-groups 
with infinitely generated centers.

It is also interesting to note that for the order
 of any C-group less than 1000 there is a 
C-group in the class described above. 
Does this hold for orders greater than 1000? That is, 
can we have a C-group of an order not covered by our class?
\\
Kazem Mahdavi\\
Department of Mathematics\\
SUNY Potsdam\\
Potsdam NY 13676\\
mahdavk@potsdam.edu\\

Avi Silterra\\
5714A Beacon St\\
Pittsburgh PA 15217\\
ajs2@andrew.cmu.edu\\

Mihai\\

Margarethe

\end{document}